\documentclass[ 11pt] {article}
\usepackage{a4,  psfrag, epsfig, epsf, amscd, amsfonts, amsmath, amssymb, amstext, amsthm, enumerate}
\begin{document}

\newtheorem{satz}{Theorem}[section]
\newtheorem{defin}[satz]{Definition}
\newtheorem{bem}[satz]{Remark}
\newtheorem{hl}[satz]{Proposition}
\newtheorem{lem}[satz]{Lemma}
\newtheorem{ko}[satz]{Corollary}

\title{Quasi-Minimal, Pseudo-Minimal Systems and Dense Orbits}
\author{Christian Pries}
\maketitle

\begin{abstract}
We give a short discussion about a weaker form of minimality (called quasi-minimality). We call a system quasi-minimal if all dense orbits form an open set. It is hard to find examples which are not already minimal. Since elliptic behaviour makes them minimal, these systems are regarded as parabolic systems. Indeed, we show that a quasi-minimal homeomorphism on a manifold is not expansive (hyperbolic). \end{abstract}

\section[]{Introduction}

\begin{defin}
A system $(X,T)$ on a topological space $X$ is said to be quasi-minimal if there exists a non empty open set $X_M \subset X$ such that  for each $x\in X_M$ the orbit $O(x)$ is dense in $X$. $X_{M,+}$ is the the subset of all $x\in X_M$ such that $O_+(x)$ is dense in $X$. Similarly we define 
$X_{M,-}$. We call $(X_M)^c$ the exceptional set.
\end{defin}

In [Gu] a system is called  quasi-minimal if $(X_M)^c$ is at most a finite set. There and in [Ka] one can find some examples. We prefer the more general definition from [Ka] where these systems are regarded as parabolic systems. Indeed, we show that there is no quasi-minimal expansive homeomorphism on a manifold. This theorem is well-known for minimal homeomorphisms (see [Ma]). Our proof here is not the same, but we follow an idea from [Ma]. Maybe Ma$\tilde{n}\acute{e}$ knew this result. Before we regard expansive homeomorphisms, we need to reproof some facts known for minimal systems.

\begin{hl}
If $(X,f)$ is quasi-minimal on a connected, perfect and locally compact metric space $X$, then $(X_M,f)$ is totally minimal.
\end{hl}

Fix a prime number $p$. We denote the elements of $\mathbb{Z}_p$ with $[q]$. For $i\in \{0,1, \ldots , p-1 \}$ and $x\in X_M $ we set
$$ A_{[i]}^{\pm}(x):=\{ y\in X_M | \quad \exists  \{n_k\}_k \in \mathbb{Z} \quad \mbox{with} \quad pn_k+i \to \pm \infty \quad\mbox{and}\quad f^{pn_k+i}(x) \to y \}.$$

\begin{lem}
\begin{enumerate}
\item If $x\in X_{m,\pm}$ then $\bigcup_i A_{[i]}^{\pm}(x)=X_M$, 
      $f(A_{[i]}^{\pm}(x))=A_{[i+1]}^{\pm}(x)$ and all 
      $A_{[i]}^{\pm}$ are open.
\item Given $x\in X_M$ and $y\in A_{[i_0]}^{+}(x)$, then for all $i$ it holds that:
      $$ A_{[i]}^{\pm}(y)  \subset A_{[i_0+i]}^{+}(x).  $$ 
\end{enumerate}
\end{lem}

Proof: 1) The first two results are trivial. From Baire's theorem we conclude that at least one $A_{[i_0]}^{\pm}(x)$ has a non empty interior. From the second fact we conclude
this holds for all these set, since $f$ is a homeomorphism. \\ 2) Choose, for $z$ in 
$A_{[i]}^{\pm}(y)$, a sequence $b_k=pm_k+i$ as in the definition, such that $f^{b_k}(y)\to z$, and a sequence $a_k=pn_k+i_0$ as in the definition such that $f^{a_k}(x)\to y$. By choosing $a_k$ growing fast, it follows that \\ $0< c_k=a_k+b_k= p(n_k+m_k)+i+i_0 \to \infty$ and $f^{c_k}(x)\to z$.
{\hfill $\Box$ \vspace{2mm}

Proof of the proposition:  We proof that $(X_M,f^p)$ is minimal for every prime number, then it follows that $(X_M,f^n)$ is minimal, 
where $n$ is a product of prime numbers. Fix $x\in X_{M,+}$. Since $X_M$ is connected, we conclude from lemma [1.3.1] that there exists  an 
$$y\in A_{[i_0]}^{+}(x) \cap A_{[i_1]}^{+}(x)$$
where $[i_1] \ne [i_0]$ and that for all $k $ and $i\in\{i_0, i_1\}$ we have 
$$A_{[k]}^{\pm}(y) \subset A_{[k+i]}^{+}(x). $$ From this and lemma [1.3.1] we conclude that 
$$ A_{[i_0]}^{+}(x)\cap A_{[i_1]}^{+}(x)   $$ must contain $A_{[0]}^{\pm}(y)$. Gottschalk showed in [G] that there is no forward or backward minimal homeomorphism on any non compact locally compact space, hence  we can suppose that $y\in X_{M,+}$ and take an $k_0$ such that $x\in A_{[k_0]}^{+}(y)$. Again we conclude that for all $n$ and $i\in\{i_0, i_1\}$ we have:
$$ A_{[n]}^{+}(x)\subset A_{[k_0+n]}^{+}(y) \subset A_{[k_0+n+i]}^{+}(x).$$ 
We suppose that $[c_0]:= [k_0+i_0]\not= 0$ (otherwise take $i_1$) and conclude,
since $p$ is a prime number, that for all $n$
$$A_{[n]}^{+}(x)=A_{[n+c_0]}^{+}(x)= ... = A_{[0]}^{+}(x),$$ which means the orbit of $x$ under $f^p$ is dense. {\hfill $\Box$ \vspace{2mm}\\ Let us generalize the theorem.
\begin{satz}
If $(X,f)$ is quasi-minimal on a locally connected, perfect and locally compact metric space $X$ then there is a finite number of connected components $\{C_i\}_{i\in I}$ and there exists a $k$ such that $f^k(C_i)=C_i$ for all $i$. In particular, $k$ is the number of connected components $\{C_i\}_{i\in I}$ and for every prime number $p>k$ we have $(X_M,f^p)$ minimal.
\end{satz}
Proof:
Let $\{C_i\}_{i\in I}$ be the connected components. Given $x\in X_{M,+}\cap C_0$, let $j(i)$ denote the unique element of $I$ such that $f^i(x)\in C_{j(i)}$.
We conclude $f^i(C_0)= C_{j(i)}$, since $f$ is an homeomorphism. The orbit of $x$ is dense, so there exist two $i_0 < i_1$ such that $f^{i_0}(C_0)= f^{i_1}(C_0)$ where the distance between $i_0$ and $i_1$ is minimal. So there can be only finitely many components and $k := i_1 - i_0$ must be the number of connected components $\{C_i\}_{i\in I}$.  The proof of proposition [1.2] works well in the case 
where $p$ is bigger than the number of connected components, since then there must exist an element
$$y\in A_{[i_0]}^{+}(x) \cap A_{[i_1]}^{+}(x).$$
for every $x\in X_{m,+}${\hfill $\Box$ \vspace{2mm}

\section[]{Topological Dynamics}

\begin{defin}
A homeomorphism $f:X\to X$ on a metric space $(X,d)$ is called expansive if there is a 
constant $e> 0$ such that for $x \ne y$ we have $\sup_{n\in \mathbb{Z}} \{ d(f^n(x), f^n(y)) \} \ge e$.
\end{defin}

Now we begin to proof the following theorem.

\begin{satz}
Suppose $(X,d)$ is a non trivial, locally connected and compact metric space and $(X,f)$ is expansive, then $f$ is not quasi-minimal.
\end{satz}

Given $\epsilon > 0$ and $x\in X$, we can define the local stable set $W^s_{\epsilon}(x)$ and the local unstable set $W^u_{\epsilon}(x)$ by
$$W^s_{\epsilon}(x)= \{y\in X \quad |  \quad  d(f^i(x), f^i(y))\le \epsilon , \forall \quad i\ge 0\},$$
$$W^u_{\epsilon}(x)= \{y\in X \quad |  \quad  d(f^i(x), f^i(y))\le \epsilon , \forall \quad i\le 0\}.$$

Here are some fundamental theorems.

\begin{satz}[Reddy]
If  $f:X\to X$ is an expansive homeomorphism on a compact metric space $(X,d)$, then there are constants $a>0$, $b>0$, $0< \gamma < 1$ and a compatible metric $D$ such that
for all $\epsilon \le a $ we have 
$$ y\in W^s_{\epsilon}(x) \rightarrow  D(f^i(x), f^i(y))\le b\gamma^i D(x,y)  \quad \forall \quad i\ge 0$$
$$ y\in W^u_{\epsilon}(x) \rightarrow  D(f^i(x), f^i(y))\le b\gamma^i D(x,y)  \quad \forall \quad i\le 0 $$
\end{satz}
Proof: See [Re] {\hfill $\Box$ \vspace{2mm}}

Let us define, for $\sigma = \{ u, s \}$ $C_{\epsilon}^{\sigma}(x)$, the connected component of $x$ in $W^{\sigma}_{\epsilon}(x)$ and $S_{\delta}(x):= \{  y\quad | \quad d(x,y)= \delta \}$. Note that $C_{\epsilon}^{\sigma}(x)\subset B(x,\epsilon) $. The following holds:
\begin{satz}[Hiriade]
If If  $f:X\to X$ is an expansive homeomorphism on a non trivial, connected, locally connected and compact metric space $(X,d)$ then for every $\epsilon > 0$
there is an $\delta (\epsilon) > 0$ such that 
$$ C_{\epsilon}^{\sigma}(x) \cap S_{\delta (\epsilon)}(x) \not= \emptyset .$$
\end{satz}
Proof: See proposition $C$ in [Hi].  {\hfill $\Box$ \vspace{2mm}}  \\ It is easy to see that in theorem [2.4] for all $\delta \le  \delta (\epsilon)$ 
$$ C_{\epsilon}^{\sigma}(x) \cap S_{\delta}(x) \not= \emptyset \quad \mbox{holds}.$$
Before showing our theorem we need two usefull lemmas.
\begin{lem}
Suppose $(X,D)$ is a non trivial, connected, locally connected and compact metric space and $(X,f)$ is an expansive quasi-minimal homeomorphism, where $D$ is the metric from theorem [2.3]. Moreover
$a,b$ and $\gamma$ are the constants from theorem [2.3]. Then  there exist constants $\epsilon > 0$, $c> 0$, $r> 0$ and a prime integer $p$ such that the following hold:
\begin{enumerate}
\item $Y:= X-B(X_M^c, 2\epsilon)$ contains an element $z$ such that $B(z,r) \subset Y$
\item $0< \epsilon < a$ and $4\epsilon < r $
\item $4c < \epsilon$, $8c < \frac{\epsilon}{2b}$ and $4c < \delta (\epsilon) $ where 
      $\delta (\epsilon)$ is taken from theorem [2.4]
\item $f^p$ is minimal on $X_M$
\item If $\Gamma$ is a connected compact subset of $W^{s}_{\epsilon}(x)$ such that $x\in \Gamma$ and $c\le Diam(\Gamma ) \le 2c$, then 
$f^{-p}(\Gamma) \subset W^{s}_{\frac{\epsilon}{2} }(f^{-p}(x)) $ and $Diam(f^{-p}(\Gamma) ) \ge 6c$.
\end{enumerate}
\end{lem}
Proof: Choose first $0< \epsilon < a$
so small  that 1 and 2 is satisfied for some $r> 0 $. Choose now a large prime integer $p$ such that 4 holds and moreover such that for $g,h \in X$ with 
$ g\in W^s_{\epsilon}(h)$ we have $D(f^p(g), f^p(h))\le \frac{D(g,h)}{12}$ (use [2.3]).
Choose now a small $c> 0$ such that 3 holds and moreover, if $Q$ is any set with
 $ Diam( Q) \le 2c$, then $D(f^i(g), f^i(h)) < \frac{\epsilon}{2} $ for all
 $-p \le i \le 0 $ and $g, h \in Q $ (use compactness). If $\Gamma$ is a connected compact subset of $W^{s}_{\epsilon}(x)$ such that $x\in \Gamma$ and $c\le Diam(\Gamma ) \le 2c$, then 
$f^{-p}(\Gamma) \subset W^{s}_{\frac{\epsilon}{2} }(f^{-p}(x)) $. This follows from theorem [2.3], since we have   by construction  $D(f^i(g), f^i(h)) < \frac{\epsilon}{2} $ for all
 $p \ge i \ge 0 $ and $g, h \in f^{-p}(\Gamma) $. Moreover we have $Diam (\Gamma)< 2c$
and therefore for all $i\ge 0$ and $l \in \Gamma$
$$ D(f^i(l), f^i(x)) \le  b\gamma^i D(l,x) < b\gamma^i 2c < \frac{\epsilon}{4}, $$
hence $D(f^i(g), f^i(h)) < \frac{\epsilon}{2} $ for all $i\ge 0$ and $g,h \in \Gamma$, thus $f^{-p}(\Gamma) \subset W^s_{ \frac{\epsilon}{2} }(f^{-p}(x)) $.
 If $Diam(f^{-p}(\Gamma) ) < 6c$, then $Diam (\Gamma ) < \frac{c}{2}$.
{\hfill $\Box$ \vspace{2mm}}

\begin{lem}
Given a non trivial, connected, locally connected and compact metric space $(X,D)$, the numbers 
$c, \epsilon $ as in lemma [2.5]  and a connected compact subset $\Gamma$ of $W^{s}_{\frac{\epsilon}{2} }(x)$ with $x\in \Gamma$ and $ Diam(\Gamma ) \ge 6c$, then there are two points $\alpha, \beta \in \Gamma$
and two compact connected sets $\alpha\in \Gamma_{\alpha}, \beta \in \Gamma_{\beta}$ such that the following holds:
\begin{enumerate}
\item $\Gamma_{\alpha}, \Gamma_{\beta} \subset \Gamma$
\item $\inf\{ d(g,h)  \quad | \quad  g\in \Gamma_{\alpha}, 
      h\in \Gamma_{\beta} \}        > \frac{c}{2}$
\item $\Gamma_{\alpha}\subset W^{s}_{\epsilon}(\alpha) $ 
      and $\Gamma_{\beta}\subset W^{s}_{\epsilon}(\beta) $ 
\item $c\le Diam(\Gamma_{\alpha, \beta} ) \le 2c$.
\end{enumerate}
\end{lem}
Proof: Since $\Gamma$ is compact and connected, choose two points $\alpha, \beta \in \Gamma$ such that $d(\alpha, \beta) = 4c$. Let $\Gamma_{\alpha} $
be the  connected component of $\alpha $ in $B(\alpha, c) \cap \Gamma$
and define $\Gamma_{\beta } $ analogously. 
It is clear that 1, 2 holds and 4 follows from $c \le \delta ( \epsilon )$. 3 follows from the triangle inequality as in the proof of [2.5]. {\hfill $\Box$ \vspace{2mm}}

Proof of the theorem: First we know that there are at most finitely many connected components,
since $X$ is locally connected and $f$ is quasi-minimal. W.l.o.g $f$ is strictly quasi-minimal, otherwise we apply [Ma]. By applying lemma [2.4] we choose
an $x$ in $B:=B(z, r)$ such that $W:=W^{s}_{\epsilon }(x)$ lies in $B$ and there is an open set $U$ with $DiamU < \frac{c}{2}$ with $W\cap U = \emptyset$.
Moreover we can assume that $O_+(x)$ is not dense hence $O_+(y)$ is not dense
for all $y\in W$ and so the backward orbit of $y$ under $f^p$ is dense for all $y\in W$,
otherwise we choose $\epsilon$ smaller. Choose by applying the previous lemmas a compact connected set $\Gamma$ in $W$
containg $x$ with $c\le Diam (\Gamma ) \le 2c$. Indeed, $C^s_{\epsilon}$ intersects 
  $  S_{\delta}(x) $ for all $\delta \le  \delta (\epsilon)$. But $4c< \delta (\epsilon)$ and therefore the connected component of $x$ (denoted by $\Gamma$) in $C^s_{\epsilon} \cap B(x,c)$ satisfies $c\le Diam (\Gamma ) \le 2c$. We know from lemma [2.4] that 
$f^{-p}(\Gamma) \subset W^{s}_{\frac{\epsilon}{2} }(f^{-p}(x)) $ and $Diam(f^{-p}(\Gamma) ) \ge 6c$, so  we apply lemma [2.5] to $f^{-p}(\Gamma)$ to  get a set $\Gamma_1$ with
$\Gamma_1 \subset f^{-p}(\Gamma)$, $x_1 \in \Gamma_1 \subset W^{s}_{\epsilon}(x_1)$
$c\le Diam(\Gamma) \le 2c $ and $\Gamma_1 \cap U = \emptyset $. Again we conclude that  
$f^{-p}(\Gamma_1) \subset W^{s}_{\frac{\epsilon}{2} }(f^{-p}(x)) $ and $Diam(f^{-p}(\Gamma_1 ) ) \ge 6c$. Repeating this construction, we get a sequence of compact connected sets $\Gamma_i$ with
$\Gamma_0 = \Gamma $, $\Gamma_i \cap U = \emptyset $ and 
$f^p(\Gamma_{i+1}) \subset \Gamma_i$. Take a point $w$ in $\bigcap_i f^{ip}(\Gamma_i) \subset W$. Then we have that  $O_{-,f^p}(w)$ does not intersect 
$U$, hence $f$ is not quasi-minimal. {\hfill $\Box$ \vspace{2mm}} \\ One could ask why the proof does not work in the case where the map is only transitive:
We need quasi-minimality to have $W:=W^{s}_{\epsilon }(x) \subset X_M$.

\section[]{Algebraic Dynamics}
We make a short excurse to algebraric dynamics. We think this is a quite natural question about quasi-minimal systems in this context.
\begin{satz}
Let $(G, \circ ) $ be a Lie group and $ \tau $ an automorphism. If $(G,\tau)$ is quasi-minimal, then $G$ is trivial.
\end{satz}

\begin{lem}
If $(G,\tau)$ is quasi-minimal, then $G$ is compact and connected.
\end{lem}
Proof: Aoki showed in [A] that an automorphism having a dense orbit on a 
 locally compact metric group $X$ implies that $X$ is compact. Therefore $G$ is compact. Since $e$ is a fixed point, $G$ is connected.
{\hfill $\Box$ \vspace{2mm}

Proof of theorem [3.1] : W.l.o.g. $G$ is compact and connected. Let us suppose that $G$ is abelian, then $G$ is a torus.
It is well-known, that for every automorphism on $G$ the periodic orbits are dense, hence $G$ is a finite set, therefore $G$ is trivial. We proof the lemma now by induction on the dimension. It holds for dim 0 and 1. Given now a non abelian Lie group of dim $G=n$. $\tau$ is ergodic by [R], since there is a dense orbit and therefore from lemma 1 in [K] we conclude that $G$ is nilpotent and the center $Z$ of $G$ has positive dimension. The induced automorphism on the lower dimensional Lie group $G/Z$ is quasi-minimal and therefore is $G=Z$, so $G$ is abelian, hence trivial.{\hfill $\Box$ 
Christian Pries \\ Fakult$\ddot{a}$t f$\ddot{u}$r Mathematik\\ Ruhr-Universit$\ddot{a}$t 
Bochum\\ Universit$\ddot{a}$tsstr. 150 \\ 44780 
Bochum \\ Germany\\ Christian.Pries@rub.de

\end{document}